\documentclass{article}
\usepackage{amssymb}
\usepackage{amsfonts}

\newtheorem{theorem}{Theorem}

\newtheorem{corollary}[theorem]{Corollary}

\newtheorem{definition}[theorem]{Definition}
\newtheorem{example}[theorem]{Example}

\newtheorem{lemma}[theorem]{Lemma}

\newtheorem{remark}[theorem]{Remark}

\newenvironment{proof}[1][Proof]{\noindent\textbf{#1.} }{\ \rule{0.5em}{0.5em}}
\input{tcilatex}
\begin{document}

\title{\textbf{On Topological Bihyperbolic Modules}}
\author{Soumen Mondal$^{1}$, Chinmay Ghosh$^{2}$, Sanjib Kumar Datta$^{3}$ \\
$^{1}$28, Dolua Dakshinpara Haridas Primary School\\
Beldanga, Murshidabad\\
Pin-742133\\
West Bengal, India\\
mondalsoumen79@gmail.com\\
$^{2}$Department of Mathematics\\
Kazi Nazrul University\\
Nazrul Road, P.O.- Kalla C.H.\\
Asansol-713340, West Bengal, India \\
chinmayarp@gmail.com \\
$^{3}$Department of Mathematics\\
University of Kalyani\\
P.O.-Kalyani, Dist-Nadia, PIN-741235,\\
West Bengal, India\\
sanjibdatta05@gmail.com}
\date{}
\maketitle

\begin{abstract}
In this paper we introduce topological modules over the ring of bihyperbolic
numbers. We discuss bihyperbolic convexity, bihyperbolic-valued seminorms
and bihyperbolic-valued Minkowski functionals in topological bihyperbolic
modules. Finally we introduce locally bihyperbolic convex modules.


\end{abstract}

\section{Introduction}

In $1882$, Corrado Segre $\cite{Se}$ introduced bicomplex numbers. Bicomplex
numbers are generalization of complex numbers by four real numbers and form
a commutative ring with divisors of zero. With the discovery of bicomplex
numbers (Tessarines), a new number system has been found which is called a
real Tessarines and defined as the set of $a+\mathbf{jc,}$ where $a,c\in 
\mathbb{R}
,$ $\mathbf{j}^{2}=1,$ $\mathbf{j\notin 
\mathbb{R}
}.$ The real Tessarine numbers are called hyperbolic numbers.

In $2002,$ S. Olariu introduced hyperbolic four complex numbers $\cite{Ol}$
which are also called bihyperbolic numbers in $\cite{Po}$. Many properties
of such numbers have been discovered during the last few year. These numbers
form a commutative ring with divisors of zero. Algebraic properties of
bihyperbolic numbers have been discussed in $\cite{Bi}$.

Topological vector spaces are one of the basic structures investigated in
functional analysis. The bicomplex version of topological module spaces was
introduced in $\cite{11}$ and \ some basic concepts and results on it have
been discussed in $\cite{16}$.

In this paper we define topological modules over the ring of bihyperbolic
numbers and discuss some basic concepts and results on it. We also discuss
bihyperbolic-valued seminorm in section \ref{S4}, bihyperbolic-valued
Minkowski functionals in section \ref{S5}\ and locally bihyperbolic convex
modules in section \ref{S6}.

\section{A Review of Bihyperbolic Numbers}

\ In this section we state some basic facts about bihyperbolic numbers.

\ \ \ \ \ Bihyperbolic numbers (also called canonical hyperbolic quaternions
or hyperbolic four complex numbers) set is defined by 
\[
H_{2}:=\{\zeta =x+y\mathbf{j}_{1}+z\mathbf{j}_{2}+w\mathbf{j}_{3}:x,y,z,w\in 
\mathbb{R}
;\mathbf{j}_{1}\mathbf{,j}_{2},\mathbf{j}_{3}\notin 
\mathbb{R}
\} 
\]

where the multiplication is given by the following rules%
\[
\mathbf{j}_{1}^{2}=\mathbf{j}_{2}^{2}=\mathbf{j}_{3}^{2}=1,\text{ }\mathbf{j}%
_{1}\mathbf{j}_{2}=\mathbf{j}_{2}\mathbf{j}_{1}=\mathbf{j}_{3},\text{ }%
\mathbf{j}_{2}\mathbf{j}_{3}=\mathbf{j}_{3}\mathbf{j}_{2}=\mathbf{j}_{1},%
\text{ }\mathbf{j}_{3}\mathbf{j}_{1}=\mathbf{j}_{1}\mathbf{j}_{3}=\mathbf{j}%
_{2}. 
\]

The canonical form of $\zeta =x+yi+zj+wk$ $\in H_{2}$ is 
\[
\zeta =\lambda _{1}(\zeta )\mathbf{e}_{1}+\lambda _{2}(\zeta )\mathbf{e}%
_{2}+\lambda _{3}(\zeta )\mathbf{e}_{3}+\lambda _{4}(\zeta )\mathbf{e}_{4} 
\]

where%
\[
\lambda _{1}(\zeta )=(x+y+z+w),\text{ }\lambda _{2}(\zeta )=(x-y+z-w),\text{ 
}\lambda _{3}(\zeta )=(x+y-z-w),\text{ }\lambda _{4}(\zeta )=(x-y-z+w) 
\]

and%
\[
\mathbf{e}_{1}=\frac{(1+\mathbf{j}_{1}+\mathbf{j}_{2}+\mathbf{j}_{3})}{4},%
\text{ }\mathbf{e}_{2}=\frac{(1-\mathbf{j}_{1}+\mathbf{j}_{2}-\mathbf{j}_{3})%
}{4},\text{ }\mathbf{e}_{3}=\frac{(1+\mathbf{j}_{1}-\mathbf{j}_{2}-\mathbf{j}%
_{3})}{4},\text{ }\mathbf{e}_{4}=\frac{(1-\mathbf{j}_{1}-\mathbf{j}_{2}+%
\mathbf{j}_{3})}{4}. 
\]

The bihyperbolic numbers $\mathbf{e}_{1},\mathbf{e}_{2},\mathbf{e}_{3},%
\mathbf{e}_{4}$ have the following properties%
\[
\mathbf{e}_{i}\mathbf{e}_{j}=\left\{ 
\begin{array}{c}
0,i\neq j\text{ } \\ 
\mathbf{e}_{i},i=j%
\end{array}%
\right. i,j\in \{1,2,3,4\} 
\]

and%
\[
\sum\limits_{i=1}^{4}\mathbf{e}_{i}=1. 
\]

The ring of bihyperbolic numbers $(H_{2},+,.)$ is a commutative ring. The
inverse of a bihyperbolic number $\zeta =x+y\mathbf{j}_{1}+z\mathbf{j}_{2}+w%
\mathbf{j}_{3}$ exist if $\lambda _{k}(\zeta )\neq 0,\forall k=1,2,3,4.$ The
set of non-invertible bihyperbolic numbers is called null cone as%
\[
NC=\{\zeta :\lambda _{1}(\zeta )\lambda _{2}(\zeta )\lambda _{3}(\zeta )%
\text{ }\lambda _{4}(\zeta )=0\}. 
\]

A non zero bihyperbolic number $\zeta =x+y\mathbf{j}_{1}+z\mathbf{j}_{2}+w%
\mathbf{j}_{3}$ is called zero divisor if there exist a non zero
bihyperbolic number $\zeta ^{\prime }=x^{\prime }+y^{\prime }\mathbf{j}%
_{1}+z^{\prime }\mathbf{j}_{2}+w^{\prime }\mathbf{j}_{3}$ such that $\zeta
\zeta ^{\prime }=0.$ Thus zero divisors exist if $\zeta ,\zeta ^{\prime }$
satisfy the following equations%
\begin{eqnarray*}
x+y+z+w &=&0\text{ and }x^{\prime }=y^{\prime }=z^{\prime }=w^{\prime } \\
&&\text{or} \\
x-y+z-w &=&0\text{ and }x^{\prime }=-y^{\prime }=z^{\prime }=-w^{\prime } \\
&&\text{or} \\
x+y-z-w &=&0\text{ and }x^{\prime }=y^{\prime }=-z^{\prime }=-w^{\prime } \\
&&\text{or} \\
x-y-z+w &=&0\text{ and }x^{\prime }=-y^{\prime }=-z^{\prime }=w^{\prime }.
\end{eqnarray*}

It is easy to verify that the set of zero divisors of the ring $H_{2}$ =$%
NC-\{0\}=NC^{\ast }$ (say), where $0=0+0\mathbf{j}_{1}+0\mathbf{j}_{2}+0%
\mathbf{j}_{3}\in H_{2}.$

It is clear that $\mathbf{e}_{i}$ $(i\in \{1,2,3,4\})$ are zero divisors of
the ring $H_{2}.$

Thus the sets%
\[
H_{2}(\mathbf{e}_{i}):=\mathbf{e}_{i}.H_{2}\text{ ,}i\in \{1,2,3,4\} 
\]

are (principal) ideals in the ring $H_{2}$ such that%
\[
H_{2}(\mathbf{e}_{i})\cap H_{2}(\mathbf{e}_{j})=\phi ,\text{ }i\neq j;i,j\in
\{1,2,3,4\}. 
\]

and $H_{2}$ can be decomposed into the following direct sum%
\[
H_{2}=H_{2}(\mathbf{e}_{1})\oplus H_{2}(\mathbf{e}_{2})\oplus H_{2}(\mathbf{e%
}_{3})\oplus H_{2}(\mathbf{e}_{4}). 
\]

Observe that%
\[
H_{2}(\mathbf{e}_{i})=\{r\mathbf{e}_{i}:r\in 
\mathbb{R}
\}=%
\mathbb{R}
\mathbf{e}_{i},\text{ }i\in \{1,2,3,4\}. 
\]

\begin{remark}
We have the following useful property

$(1)$ $\zeta \in H_{2}(\mathbf{e}_{i})$ if and only if $\zeta \mathbf{e}%
_{i}=\zeta $, for $i\in \{1,2,3,4\}.$

$(2)$ $\zeta \in H_{2}(\mathbf{e}_{i})\oplus H_{2}(\mathbf{e}_{j})$ if and
only if $\zeta (\mathbf{e}_{i}+\mathbf{e}_{j})=\zeta ,$ for $i\neq j;$ $%
i,j\in \{1,2,3,4\}.$

$(3)$ $\zeta \in H_{2}(\mathbf{e}_{i})\oplus H_{2}(\mathbf{e}_{j})\oplus
H_{2}(\mathbf{e}_{k})$ if and only if $\zeta (\mathbf{e}_{i}+\mathbf{e}_{j}+%
\mathbf{e}_{k})=\zeta ,$ for $i\neq j\neq k;$ $i,j,k\in \{1,2,3,4\}.$
\end{remark}

The set of non-negative bihyperbolic number is%
\[
H_{2}^{+}:=\{\zeta =w_{1}\mathbf{e}_{1}+w_{2}\mathbf{e}_{2}+w_{3}\mathbf{e}%
_{3}+w_{4}\mathbf{e}_{4}:w_{1},w_{2},w_{3},w_{4}\geq 0\}. 
\]

The binary relation on $H_{2}$ defined by:%
\[
\zeta \preceq \varphi \text{ if and only if }\lambda _{k}(\zeta )\leq
\lambda _{k}(\varphi )\text{ for all }k\in \{1,2,3,4\} 
\]

is a partial order on $H_{2}.$ If we take $\alpha ,\beta \in 
\mathbb{R}
,$ then $\alpha \preceq \beta $ if and only if $\alpha \leq \beta .$

For $\zeta ,\eta ,\xi ,\vartheta \in H_{2},$ it is easy to verify that

(1) If $\zeta \preceq \eta $ and $\xi \in H_{2}^{+},$ then $\zeta \xi
\preceq \eta \xi .$

(2) If $\zeta \preceq \eta $ and $\xi \preceq \vartheta ,$ then $\zeta +\xi
\preceq \eta +\vartheta .$

(3) If $\zeta \preceq \eta ,$ then $-\eta \preceq -\zeta .$

The bihyperbolic-valued modulus of a number $\zeta $ is defined by 
\[
\left\vert \zeta \right\vert =\left\vert \lambda _{1}(\zeta )\right\vert 
\mathbf{e}_{1}+\left\vert \lambda _{2}(\zeta )\right\vert \mathbf{e}%
_{2}+\left\vert \lambda _{3}(\zeta )\right\vert \mathbf{e}_{3}+\left\vert
\lambda _{4}(\zeta )\right\vert \mathbf{e}_{4}. 
\]

This satisfies the following properties:

(1) $\left\vert \zeta \right\vert =0$ if and only if $\zeta =0.$

(2) $\left\vert \zeta \varphi \right\vert =\left\vert \zeta \right\vert
.\left\vert \varphi \right\vert .$

(3) $\left\vert \zeta +\varphi \right\vert \preceq \left\vert \zeta
\right\vert +\left\vert \varphi \right\vert $ for any $\zeta ,\varphi \in
H_{2}.$

Let $A\subset H_{2},$ if there exists $M\in H_{2}^{+}$ such that $\left\vert
x\right\vert \preceq M$ \ \ $\forall x\in A,$ we say that $A$ is $H_{2}-$%
bounded set.

If $A\subset H_{2}$ is $H_{2}-$bounded from above, then the $H_{2}-$supremum
of $A$ is defined as%
\[
\sup_{H_{2}}A=\sum\limits_{i=1}^{4}\sup A_{i}\mathbf{e}_{i}, 
\]

where $A_{i}=\{a_{i}:\sum\limits_{i=1}^{4}a_{i}\mathbf{e}_{i}\in A\}.$

Similarly $H_{2}-$infimum of a $H_{2}-$bounded below set $A$ is defined as%
\[
\inf_{H_{2}}A=\sum\limits_{i=1}^{4}\inf A_{i}\mathbf{e}_{i}, 
\]

where $A_{i}$ are defined as above.

\section{Topological Bihyperbolic Modules}

Topological bicomplex modules have been introduced in $\cite{11}$. In this
section, we introduce topological bihyperbolic modules, the concept of
balancedness, convexity and absorbedness in bihyperbolic modules and discuss
some of their properties.

\begin{definition}
Let $X$ be a $H_{2}$ module and $\tau $ be a Hausdorff topology on $X$ such
that the operations

$(i)$ $+:X\times X\longrightarrow X$ and

$(ii)\,$\ $.:H_{2}\times X\longrightarrow X$

are continuous. Then the pair $(X,\tau )$ is called a topological
bihyperbolic module or topological $H_{2}-$ module.
\end{definition}

\begin{example}
Every $H_{2}-$module with $H_{2}-$valued norm is a topological $H_{2}-$%
module.
\end{example}

\begin{remark}
Let $(X,\tau )$ be a topological $H_{2}-$ module. Write%
\[
X=\sum\limits_{i=1}^{4}\mathbf{e}_{i}X 
\]

where $X_{i}=\mathbf{e}_{i}X$ are $%
\mathbb{R}
-$vector spaces. Then $\tau _{l}=\{\mathbf{e}_{l}G:G\in \tau \}$ is a
Hausdorff topology on $X_{l}$ and so $(X_{l},\tau _{l})$ is a topological $%
\mathbb{R}
-$vector space for $l\in \{1,2,3,4\}.$
\end{remark}

\begin{definition}
Let $X$ be a $H_{2}$ module and a function $\left\Vert .\right\Vert
_{H_{2}}:X\longrightarrow H_{2}^{+}$ such that

$1.\left\Vert x\right\Vert _{H_{2}}=0\Leftrightarrow x=0,$

$2.\left\Vert \lambda .x\right\Vert _{H_{2}}=\left\vert \lambda \right\vert
\left\Vert x\right\Vert _{H_{2}}$ $\forall x\in X,$ $\forall \lambda \in
H_{2},$

$3.\left\Vert x+y\right\Vert _{H_{2}}\preceq \left\Vert x\right\Vert
_{H_{2}}+\left\Vert y\right\Vert _{H_{2}},$

is called bihyperbolic-valued or $H_{2}-$valued norm on $X.$
\end{definition}

If $\left\Vert .\right\Vert _{i}$ are $%
\mathbb{R}
-$valued norm on $X_{i}$ for $i\in \{1,2,3,4\},$ then $X$ can be endowed
canonically with $H_{2}-$valued norm given by the formula 
\[
\left\Vert x\right\Vert _{H_{2}}=\left\Vert \sum\limits_{i=1}^{4}\mathbf{e}%
_{i}x_{i}\right\Vert _{H_{2}}=\sum\limits_{i=1}^{4}\left\Vert
x_{i}\right\Vert _{i}\mathbf{e}_{i}. 
\]

\begin{lemma}
\label{L1}Let $X$ be a topological $H_{2}-$module. Then, for any $y\in X,$
the map $T_{y}:X\longrightarrow X$ defined by 
\[
T_{y}(x)=x+y\text{ \ for each }x\in X, 
\]

is a homeomorphism.

\begin{proof}
The proof is similar to ($\cite{16},$ Lemma $2.5$).
\end{proof}
\end{lemma}

\begin{lemma}
\label{L2}Let $X$ be a topological $H_{2}-$module. Then for any $\lambda \in
H_{2}\backslash NC,$ the map $M_{\lambda }:X\longrightarrow X$ defined by 
\[
M_{\lambda }(x)=\lambda .x\text{ \ \ for each }x\in X, 
\]

is a homeomorphism.

\begin{proof}
The proof is similar to ($\cite{16}$, Lemma $2.6$).
\end{proof}
\end{lemma}

\begin{definition}
Let $B$ be a subset of a $H_{2}-$module $X.$ Then $B$ is called a $H_{2}-$%
balanced set if for any $x\in B$ and $\lambda \in H_{2}$ with $\left\vert
\lambda \right\vert \preceq 1,$ $\lambda x\in B.$
\end{definition}

In other words, $\lambda B\subseteq B$ for any $\lambda \in H_{2},$ $%
\left\vert \lambda \right\vert \preceq 1.$ It is obvious that if$\ B$ is $%
H_{2}-$balanced set, then $0\in H_{2}.$

\begin{theorem}
\label{T1}Let $B$ be a $H_{2}-$balanced subset of a $H_{2}-$module $X.$ Then

$(i)$ $\lambda B=B$ whenever $\lambda \in H_{2}$ with $\left\vert \lambda
\right\vert =1.$

$(ii)$ $\lambda B=\left\vert \lambda \right\vert B$ for each $\lambda \in
H_{2}\backslash NC^{\ast }.$

\begin{proof}
$(i)$ Let $\lambda \in H_{2}$ with $\left\vert \lambda \right\vert =1.$
Since $B$ is $H_{2}-$balanced, $\lambda B\subseteq B.$ Writing $\lambda
=\sum\limits_{i=1}^{4}\lambda _{i}\mathbf{e}_{i},$ we have $\left\vert
\lambda _{i}\right\vert =1$ for each $i\in \{1,2,3,4\}.$

Therefore 
\[
\left\vert \frac{1}{\lambda }\right\vert =\frac{1}{\left\vert \lambda
\right\vert }=\left\vert \lambda \right\vert
^{-1}=\sum\limits_{i=1}^{4}\left\vert \lambda _{i}\right\vert ^{-1}\mathbf{e}%
_{i}=1. 
\]

So,%
\[
\frac{1}{\lambda }B\subseteq B\Longrightarrow B\subseteq \lambda B. 
\]

$(ii)$ Let $\lambda \in H_{2}\backslash NC^{\ast }.$ If $\lambda =0,$ then $%
\lambda B=\left\vert \lambda \right\vert B.$ Now let $\lambda \neq 0.$
Writing $\lambda =\sum\limits_{i=1}^{4}\lambda _{i}\mathbf{e}_{i},$ we obtain%
\[
\frac{\lambda }{\left\vert \lambda \right\vert }=\sum\limits_{i=1}^{4}\frac{%
\lambda _{i}}{\left\vert \lambda _{i}\right\vert }\mathbf{e}_{i}. 
\]

Hence%
\[
\left\vert \frac{\lambda }{\left\vert \lambda \right\vert }\right\vert =1. 
\]

So by $(i),$ we have%
\[
\frac{\lambda }{\left\vert \lambda \right\vert }B=B\Longrightarrow \lambda
B=\left\vert \lambda \right\vert B. 
\]
\end{proof}
\end{theorem}

\begin{theorem}
\label{T2}Let $B$ be a $H_{2}-$balanced subset of $H_{2}-$module $X.$ Then

$(i)$ $\mathbf{e}_{i}B$ are balanced sets in $%
\mathbb{R}
-$vector spaces $\mathbf{e}_{i}X$,

$(ii)$ $\mathbf{e}_{i}B\subset B$ for all $i\in \{1,2,3,4\}.$

\begin{proof}
$(i)$ Let $x\in \mathbf{e}_{i}B$ and $a\in 
\mathbb{R}
$ such that $\left\vert a\right\vert \leq 1.$ Then there exists $x^{\prime }$
and $a^{\prime }\in H_{2}$ with $\left\vert a^{\prime }\right\vert \preceq 1$
such that $x=\mathbf{e}_{i}x^{\prime }$\ and $a=\mathbf{e}_{i}a^{\prime }.$

Since $B$ is $H_{2}-$balanced, $a^{\prime }x^{\prime }$ $\in B.$ Thus $ax=a%
\mathbf{e}_{i}x^{\prime }=\mathbf{e}_{i}a^{\prime }$ $x^{\prime }\in \mathbf{%
e}_{i}B,$ showing that $\mathbf{e}_{i}B$ is balanced set in $%
\mathbb{R}
-$vector space $\mathbf{e}_{i}X.$

$(ii)$ Let $x\in \mathbf{e}_{i}B.$ Then there is an $x^{\prime }\in B$ such
that $x=$ $\mathbf{e}_{i}x^{\prime }.$ Since $x^{\prime }\in B,$ by $H_{2}-$%
balancedness of $B,$ $\lambda x^{\prime }\in B$ for any $\lambda \in H_{2}$
with $\left\vert \lambda \right\vert \preceq 1.$

In particular, taking $\lambda =\mathbf{e}_{i},$ we get $x=\mathbf{e}%
_{i}x^{\prime }\in B.$ Thus $\mathbf{e}_{i}B\subset B.$
\end{proof}
\end{theorem}

\begin{definition}
Let $B$ be a subset of a $H_{2}-$module $X.$ Then $B$ is called a $H_{2}-$%
convex set if $x,y\in X$ and $\lambda \in H_{2}^{+}$ satisfying $0\preceq
\lambda \preceq 1$ implies that $\lambda x+(1-\lambda )y\in B.$
\end{definition}

\begin{definition}
Let $\{x_{1},x_{2},...,x_{n}\}$ be a subset of a $H_{2}-$module $X.$ Then
the linear combinations $\sum\limits_{i=1}^{n}a_{i}x_{i}$ in which $a_{i}\in
H_{2}^{+}$ and $\sum\limits_{i=1}^{n}a_{i}=1$ are called $H_{2}-$convex
combinations of the $x_{i}$'s.
\end{definition}

In this terminology, a $H_{2}-$convex is one that contains all its $H_{2}-$%
convex combinations.

\begin{theorem}
Let $B$ be a $H_{2}-$convex subset of a $H_{2}-$module $X.$ Then

$(i)$ $\mathbf{e}_{i}B$ are convex sets in $%
\mathbb{R}
-$vector spaces $\mathbf{e}_{i}X,$ for $i\in \{1,2,3,4\}.$

$(ii)$ $\mathbf{e}_{i}B\subset B$ whenever $0\in B,$ for all $i\in
\{1,2,3,4\}.$

\begin{proof}
$(i)$ Take $x,y\in \mathbf{e}_{i}B$, then there exist $x^{\prime },y^{\prime
}\in B$ such that $x=\mathbf{e}_{i}x^{\prime },$ $y=\mathbf{e}_{i}y^{\prime
}.$ Take $\lambda _{i}\in \lbrack 0,1]$ be such that $\lambda
=\sum\limits_{i=1}^{4}\lambda _{i}\mathbf{e}_{i}$ satisfy $0\preceq \lambda
\preceq 1.$ Since $B$ is $H_{2}-$convex, $\lambda x^{\prime }+(1-\lambda
)y^{\prime }\in B.$ Hence $\mathbf{e}_{i}\left( \lambda x^{\prime
}+(1-\lambda )y^{\prime }\right) =\lambda _{i}x+(1-\lambda _{i})y\in \mathbf{%
e}_{i}B,$ proving that $\mathbf{e}_{i}B$ are convex.

$(ii)$ Given $x\in \mathbf{e}_{i}B,$ take as before $x^{\prime }\in B$ such
that $x=\mathbf{e}_{i}x^{\prime }.$ Since $0\in B,$ for any $\lambda
=\sum\limits_{i=1}^{4}\lambda _{i}\mathbf{e}_{i}$ $\in H_{2}$ with $0\preceq
\lambda \preceq 1,$ there follows that $\lambda x^{\prime }\in B.$

In particular, taking $\lambda =\mathbf{e}_{i},$ one has $\mathbf{e}%
_{i}x^{\prime }=x\in B,$ i.e., $\mathbf{e}_{i}B\subset B.$
\end{proof}
\end{theorem}

The following lemma is easy to prove:

\begin{lemma}
In a $H_{2}-$module $X,$ if $\{B_{l}:l\in \Delta \}$ is a collection of $%
H_{2}-$convex sets, then $\cap _{l}B_{l}$ is $H_{2}-$convex.
\end{lemma}

\begin{theorem}
\label{T4}Let $B$ be a $H_{2}-$convex subset of $H_{2}-$module $X.$ Then $B$
can be written as $B=\sum\limits_{i=1}^{4}\mathbf{e}_{i}B.$

\begin{proof}
Let $x\in B,$ then $\mathbf{e}_{i}x\in \mathbf{e}_{i}B$ $\forall i\in
\{1,2,3,4\}.$

Therefore 
\[
x=\left( \sum\limits_{i=1}^{4}\mathbf{e}_{i}\right) x=\sum\limits_{i=1}^{4}%
\mathbf{e}_{i}x\in \sum\limits_{i=1}^{4}\mathbf{e}_{i}B. 
\]

Hence $B\subset \sum\limits_{i=1}^{4}\mathbf{e}_{i}B.$

Now let \thinspace $x_{i}\in \mathbf{e}_{i}B.$ Then there exist $%
x_{i}^{\prime }\in B$ such that $x=\mathbf{e}_{i}x_{i}^{\prime }.$

Now 
\[
\sum\limits_{i=1}^{4}x_{i}=\sum\limits_{i=1}^{4}\mathbf{e}_{i}x_{i}^{\prime
} 
\]

and it is a $H_{2}-$convex combination of elements of $B.$

Since $B$ is $H_{2}-$convex, then $\sum\limits_{i=1}^{4}x_{i}\in B.$ Thus $%
\sum\limits_{i=1}^{4}\mathbf{e}_{i}B\subset B.$
\end{proof}
\end{theorem}

\begin{theorem}
\label{T5}Let $X$ be a $H_{2}-$module and $B\subset X.$ If $\mathbf{e}_{i}B$
are convex sets in $%
\mathbb{R}
-$vector spaces $\mathbf{e}_{i}X,$ for each $i\in \{1,2,3,4\},$ then $%
\sum\limits_{i=1}^{4}\mathbf{e}_{i}B$ is a $H_{2}-$convex subset of $X.$

\begin{proof}
Let $x,y\in \sum\limits_{i=1}^{4}\mathbf{e}_{i}B$ and $0\preceq \lambda
\preceq 1.$ Write 
\[
x=\sum\limits_{i=1}^{4}\mathbf{e}_{i}x_{i},\text{ }y=\sum\limits_{i=1}^{4}%
\mathbf{e}_{i}y_{i}\text{ and }\lambda =\sum\limits_{i=1}^{4}\mathbf{e}%
_{i}\lambda _{i}, 
\]

where $\mathbf{e}_{i}x_{i},$ $\mathbf{e}_{i}y_{i}\in \mathbf{e}_{i}B$ and $%
0\leq \lambda _{i}\leq 1.$

Since $\mathbf{e}_{i}B$ are convex in $%
\mathbb{R}
-$vector spaces $\mathbf{e}_{i}X,$ then we have%
\[
\mathbf{e}_{i}\lambda _{i}x_{i}+\mathbf{e}_{i}(1-\lambda _{i})y_{i}\in 
\mathbf{e}_{i}B. 
\]

Then a simple calculation shows that 
\[
\lambda x+(1-\lambda )y\in \sum\limits_{i=1}^{4}\mathbf{e}_{i}B, 
\]

showing that $\sum\limits_{i=1}^{4}\mathbf{e}_{i}B$ is $H_{2}-$convex.
\end{proof}
\end{theorem}

If $\mathbf{e}_{i}B$ are convex sets in $%
\mathbb{R}
-$vector spaces $\mathbf{e}_{i}X$ for $i\in \{1,2,3,4\},$ then $%
B=\sum\limits_{i=1}^{4}\mathbf{e}_{i}B$ may not hold:

\begin{example}
Let $X=H_{2}$ and $B=\{x=\sum\limits_{i=1}^{4}\mathbf{e}_{i}x_{i}:x_{i}\in 
\mathbb{R}
,$ $\sum\limits_{i=1}^{4}\left\vert x_{i}\right\vert <2\}.$ Then $\mathbf{e}%
_{i}B=\{\mathbf{e}_{i}x_{i}:\left\vert x_{i}\right\vert <2\}$ are convex
sets in $%
\mathbb{R}
-$vector spaces \thinspace $\mathbf{e}_{i}X.$ Now $\mathbf{e}_{i}\frac{3}{4}%
\in \mathbf{e}_{i}B,$ but $\frac{3}{4}=\sum\limits_{i=1}^{4}\mathbf{e}_{i}%
\frac{3}{4}\notin B.$ Therefore $B\neq \sum\limits_{i=1}^{4}\mathbf{e}_{i}B$
and hence $B$ is not $H_{2}-$convex.
\end{example}

\begin{theorem}
\label{T6}Let $X$ be a topological $H_{2}-$module and $B\subset X.$ Then the
following statements hold:

$(i)$ $\left( \mathbf{e}_{i}B\right) ^{\circ }=\mathbf{e}_{i}B^{\circ }$ $\
\forall i\in \{1,2,3,4\}.$

$(ii)$ $\overline{\mathbf{e}_{i}B}=\mathbf{e}_{i}\overline{B}$ \ $\forall
i\in \{1,2,3,4\}.$

\begin{proof}
The proof is similar to ($\cite{16},$ Theorem $2.16$).
\end{proof}
\end{theorem}

\begin{theorem}
Let $B$ be a $H_{2}-$convex set in a topological $H_{2}-$module $X.$ Then
the following statements hold:

$(i)$ $B^{\circ }=\sum\limits_{i=1}^{4}\mathbf{e}_{i}B^{\circ }$ and $%
\overline{B}=\sum\limits_{i=1}^{4}\mathbf{e}_{i}\overline{B}.$

$(ii)$ $B^{\circ }$ and $\overline{B}$ are $H_{2}-$convex sets.

\begin{proof}
$(i)$ Since $B$ is $H_{2}-$convex, we can write $B$ as%
\[
B=\sum\limits_{i=1}^{4}\mathbf{e}_{i}B. 
\]

Clearly $B^{\circ }\subset \sum\limits_{i=1}^{4}\mathbf{e}_{i}B^{\circ }.$

Now $\sum\limits_{i=1}^{4}\mathbf{e}_{i}B^{\circ }$ is an open set in $X$
such that $\sum\limits_{i=1}^{4}\mathbf{e}_{i}B^{\circ }\subset
\sum\limits_{i=1}^{4}\mathbf{e}_{i}B=B.$

But $B^{\circ }$ is the largest open set contained in $B.$ Therefore $%
\sum\limits_{i=1}^{4}\mathbf{e}_{i}B^{\circ }\subset B^{\circ }.$ Thus $%
B^{\circ }=\sum\limits_{i=1}^{4}\mathbf{e}_{i}B^{\circ }.$

Again, trivially%
\[
\overline{B}\subset \sum\limits_{i=1}^{4}\mathbf{e}_{i}\overline{B}. 
\]

From Theorem \ref{T6} and $\cite{Ru}$, it follows that 
\[
\sum\limits_{i=1}^{4}\mathbf{e}_{i}\overline{B}=\sum\limits_{i=1}^{4}%
\overline{\mathbf{e}_{i}B}\subset \overline{\sum\limits_{i=1}^{4}\mathbf{e}%
_{i}B}=\overline{B}. 
\]

Thus we have $\overline{B}=\sum\limits_{i=1}^{4}\mathbf{e}_{i}\overline{B}.$

$(ii)$ Since $B$ is $H_{2}-$convex, $\mathbf{e}_{i}B$ are convex sets in $%
\mathbb{R}
-$vector spaces $\mathbf{e}_{i}X$ for $i\in \{1,2,3,4\}.$ Then it follows
from $\cite{Ru},$ $(\mathbf{e}_{i}B)^{\circ }$ are convex in $%
\mathbb{R}
-$vector spaces $\mathbf{e}_{i}X$ and hence by Theorem \ref{T6}$(i),$ $%
\mathbf{e}_{i}B^{\circ }$ are convex in $%
\mathbb{R}
-$vector spaces $\mathbf{e}_{i}X.$ Now from Theorem \ref{T5} we see that $%
\sum\limits_{i=1}^{4}\mathbf{e}_{i}B^{\circ }$ is $H_{2}-$convex and by $%
(i), $ it follows that $B^{\circ }$ is $H_{2}-$convex subset of $X.$

Similarly we can prove that $\overline{B}$ is $H_{2}-$convex .
\end{proof}
\end{theorem}

\begin{theorem}
\label{T8}Let $B$ be a $H_{2}-$convex set in a topological $H_{2}-$module $X$
such that $0\in B.$ Then

$(i)$ $\mathbf{e}_{i}B+\mathbf{e}_{j}B\subset B$ for all $i\neq j,$ $i,j\in
\{1,2,3,4\}.$

$(ii)$ $\mathbf{e}_{i}B+\mathbf{e}_{j}B+\mathbf{e}_{k}B\subset B$ for all $%
i\neq j\neq k,$ \thinspace $i,j,k\in \{1,2,3,4\}.$

\begin{proof}
$(i)$ Let $x\in \mathbf{e}_{i}B+\mathbf{e}_{j}B.$ Then $x=\mathbf{e}%
_{i}x_{i}+\mathbf{e}_{j}x_{j}$ for some $x_{i},x_{j}\in B.$

Since, $B$ is $H_{2}-$convex subset, then $B=\sum\limits_{q=1}^{4}\mathbf{e}%
_{q}B.$ Also since $0\in B,$ then $0\in \mathbf{e}_{l}B$ for $l\neq i\neq j,$
$l\in \{1,2,3,4\}.$ Hence%
\[
x=\mathbf{e}_{i}x_{i}+\mathbf{e}_{j}x_{j}+0+0\in \sum\limits_{q=1}^{4}%
\mathbf{e}_{q}B=B. 
\]

So, $\mathbf{e}_{i}B+\mathbf{e}_{j}B\subset B$.

$(ii)$ Let $x\in \mathbf{e}_{i}B+\mathbf{e}_{j}B+\mathbf{e}_{k}B.$ Then $x=%
\mathbf{e}_{i}x_{i}+\mathbf{e}_{j}x_{j}+\mathbf{e}_{k}x_{k}$ for some $%
x_{i},x_{j},x_{k}\in B.$

Since, $B$ is $H_{2}-$convex subset, then $B=\sum\limits_{q=1}^{4}\mathbf{e}%
_{q}B.$ Also since $0\in B,$ then $0\in \mathbf{e}_{l}B$ for $l\neq i\neq
j\neq k,$ $l\in \{1,2,3,4\}.$ Hence%
\[
x=\mathbf{e}_{i}x_{i}+\mathbf{e}_{j}x_{j}+\mathbf{e}_{k}x_{k}+0\in
\sum\limits_{q=1}^{4}\mathbf{e}_{q}B=B. 
\]

So, $\mathbf{e}_{i}B+\mathbf{e}_{j}B+\mathbf{e}_{k}B\subset B$ .
\end{proof}
\end{theorem}

\begin{theorem}
Let $B$ be a $H_{2}-$balanced and $H_{2}-$convex set in a topological $%
H_{2}- $module $X.$ Then $\overline{B}$ is $H_{2}-$balanced and so is $%
B^{\circ }$ if $0\in B^{\circ }.$

\begin{proof}
Let $\lambda \in H_{2}$ such that $\left\vert \lambda \right\vert \preceq 1.$
If $\lambda =0,$ then $\lambda \overline{B}=\{0\}\subset \overline{B}.$ If $%
\lambda \notin NC,$ then by Lemma \ref{L2}, we have $\lambda \overline{B}=%
\overline{\lambda B}\subset \overline{B}.$

If $\lambda \in NC^{\ast }$ such that $\lambda =\lambda _{i}\mathbf{e}_{i}$
for $i\in \{1,2,3,4\},$ then $0<\lambda _{i}\leq 1.$ Then using respectively
Theorem \ref{T6} ($ii)$, ($\cite{La}$, Theorem $2.1.2$)$,$ balancedness of $%
\mathbf{e}_{i}B$ and Theorem \ref{T2} $(ii),$ we obtain%
\[
\lambda \overline{B}=\lambda _{i}\mathbf{e}_{i}\overline{B}=\lambda _{i}%
\overline{\mathbf{e}_{i}B}=\overline{\lambda _{i}\mathbf{e}_{i}B}\subset 
\overline{\mathbf{e}_{i}B}\subset \overline{B}. 
\]

If $\lambda \in NC^{\ast }$ such that $\lambda =\lambda _{i}\mathbf{e}%
_{i}+\lambda _{j}\mathbf{e}_{j}$ for $i\neq j,$ $i,j\in \{1,2,3,4\}.$ Then $%
0<\lambda _{i}\leq 1$ and $0<\lambda _{j}\leq 1.$ Then using Theorem \ref{T6}%
, ($\cite{La}$, Theorem $2.1.2$)$,$ balancedness of $\mathbf{e}_{i}B$ ,
Theorem \ref{T2} $(ii),$ ($\cite{Na},$ Theorem $4.4.1(b)$) and Theorem \ref%
{T8} $(i)$, we obtain%
\begin{eqnarray*}
\lambda \overline{B} &=&(\lambda _{i}\mathbf{e}_{i}+\lambda _{j}\mathbf{e}%
_{j})\overline{B}\subset \lambda _{i}\mathbf{e}_{i}\overline{B}+\lambda _{j}%
\mathbf{e}_{j}\overline{B} \\
&=&\lambda _{i}\overline{\mathbf{e}_{i}B}+\lambda _{j}\overline{\mathbf{e}%
_{j}B} \\
&=&\overline{\lambda _{i}\mathbf{e}_{i}B}+\overline{\lambda _{j}\mathbf{e}%
_{j}B} \\
&\subset &\overline{\mathbf{e}_{i}B}+\overline{\mathbf{e}_{j}B} \\
&\subset &\overline{\mathbf{e}_{i}B+\mathbf{e}_{j}B}\subset \overline{B}.
\end{eqnarray*}

If $\lambda \in NC^{\ast }$ such that $\lambda =\lambda _{i}\mathbf{e}%
_{i}+\lambda _{j}\mathbf{e}_{j}+\lambda _{k}\mathbf{e}_{k}$ for $i\neq j\neq
k,$ $i,j,k\in \{1,2,3,4\}.$ Then $0<\lambda _{i},\lambda _{j},\lambda
_{k}\leq 1.$ Then using Theorem \ref{T6}, ($\cite{La}$, Theorem $2.1.2$)$,$
balancedness of $\mathbf{e}_{i}B$ , Theorem \ref{T2} $(ii),$ ($\cite{Na},$
Theorem $4.4.1(b)$] and Theorem \ref{T8} $(ii)$, we obtain%
\begin{eqnarray*}
\lambda \overline{B} &=&(\lambda _{i}\mathbf{e}_{i}+\lambda _{j}\mathbf{e}%
_{j}+\lambda _{k}\mathbf{e}_{k})\overline{B} \\
&\subset &\lambda _{i}\mathbf{e}_{i}\overline{B}+\lambda _{j}\mathbf{e}_{j}%
\overline{B}+\lambda _{k}\mathbf{e}_{k}\overline{B} \\
&=&\lambda _{i}\overline{\mathbf{e}_{i}B}+\lambda _{j}\overline{\mathbf{e}%
_{j}B}+\lambda _{k}\overline{\mathbf{e}_{k}B} \\
&=&\overline{\lambda _{i}\mathbf{e}_{i}B}+\overline{\lambda _{j}\mathbf{e}%
_{j}B}+\overline{\lambda _{k}\mathbf{e}_{k}B} \\
&\subset &\overline{\mathbf{e}_{i}B}+\overline{\mathbf{e}_{j}B}+\overline{%
\mathbf{e}_{k}B} \\
&\subset &\overline{\mathbf{e}_{i}B+\mathbf{e}_{j}B+\mathbf{e}_{k}B}\subset 
\overline{B}.
\end{eqnarray*}

Hence $\overline{B}$ is $H_{2}-$balanced.

Now suppose that $0\in B^{\circ }$ and $\lambda \in H_{2}$ with $\left\vert
\lambda \right\vert \preceq 1.$ If $\lambda =0,$ then $\lambda B^{\circ
}=\{0\}\in B^{\circ }.$ If $\lambda \notin NC^{\ast },$ then by Lemma \ref%
{L2}, we have $\lambda B^{\circ }=(\lambda B)^{\circ }\subset B^{\circ }.$

If $\lambda \in NC^{\ast }$ such that $\lambda =\lambda _{i}\mathbf{e}_{i}$
for $i\in \{1,2,3,4\},$ then $0<\lambda _{i}\leq 1.$ Then using respectively
Theorem \ref{T6} ($i)$, ($\cite{La}$, Theorem $2.1.2$)$,$ balancedness of $%
\mathbf{e}_{i}B$ and Theorem \ref{T2} $(ii),$ we obtain%
\[
\lambda B^{\circ }=\lambda _{i}\mathbf{e}_{i}B^{\circ }=\lambda _{i}(\mathbf{%
e}_{i}B)^{\circ }=(\lambda _{i}\mathbf{e}_{i}B)^{\circ }\subset (\mathbf{e}%
_{i}B)^{\circ }\subset B^{\circ }. 
\]

If $\lambda \in NC^{\ast }$ such that $\lambda =\lambda _{i}\mathbf{e}%
_{i}+\lambda _{j}\mathbf{e}_{j}$ for $i\neq j,$ $i,j\in \{1,2,3,4\}.$ Then $%
0<\lambda _{i}\leq 1$ and $0<\lambda _{j}\leq 1.$ Then using Theorem \ref{T6}%
, ($\cite{La}$, Theorem $2.1.2$)$,$ balancedness of $\mathbf{e}_{i}B$ ,
Theorem \ref{T2} $(ii),$ ($\cite{Na},$ Theorem $4.4.1(e)$] and Theorem \ref%
{T8} $(i)$, we obtain

\begin{eqnarray*}
\lambda B^{\circ } &=&(\lambda _{i}\mathbf{e}_{i}+\lambda _{j}\mathbf{e}%
_{j})B^{\circ } \\
&\subset &\lambda _{i}\mathbf{e}_{i}B^{\circ }+\lambda _{j}\mathbf{e}%
_{j}B^{\circ } \\
&=&\lambda _{i}(\mathbf{e}_{i}B)^{\circ }+\lambda _{j}(\mathbf{e}%
_{i}B)^{\circ } \\
&=&(\lambda _{i}\mathbf{e}_{i}B)^{\circ }+(\lambda _{j}\mathbf{e}%
_{j}B)^{\circ } \\
&\subset &(\mathbf{e}_{i}B)^{\circ }+(\mathbf{e}_{j}B)^{\circ } \\
&\subset &(\mathbf{e}_{i}B+\mathbf{e}_{j}B)^{\circ }\subset B^{\circ }.
\end{eqnarray*}

If $\lambda \in NC^{\ast }$ such that $\lambda =\lambda _{i}\mathbf{e}%
_{i}+\lambda _{j}\mathbf{e}_{j}+\lambda _{k}\mathbf{e}_{k}$ for $i\neq j\neq
k,$ $i,j,k\in \{1,2,3,4\}.$ Then $0<\lambda _{i},\lambda _{j},\lambda
_{k}\leq 1.$ Then using Theorem \ref{T6}, ($\cite{La}$, Theorem $2.1.2$)$,$
balancedness of $\mathbf{e}_{i}B$ , Theorem \ref{T2} $(ii),$ ($\cite{Na},$
Theorem $4.4.1(e)$] and Theorem \ref{T8} $(ii)$, we obtain

\begin{eqnarray*}
\lambda B^{\circ } &=&(\lambda _{i}\mathbf{e}_{i}+\lambda _{j}\mathbf{e}%
_{j}+\lambda _{k}\mathbf{e}_{k})B^{\circ } \\
&\subset &\lambda _{i}\mathbf{e}_{i}B^{\circ }+\lambda _{j}\mathbf{e}%
_{j}B^{\circ }+\lambda _{k}\mathbf{e}_{k}B^{\circ } \\
&=&\lambda _{i}(\mathbf{e}_{i}B)^{\circ }+\lambda _{j}(\mathbf{e}%
_{i}B)^{\circ }+\lambda _{k}(\mathbf{e}_{k}B)^{\circ } \\
&=&(\lambda _{i}\mathbf{e}_{i}B)^{\circ }+(\lambda _{j}\mathbf{e}%
_{j}B)^{\circ }+(\lambda _{k}\mathbf{e}_{k}B)^{\circ } \\
&\subset &(\mathbf{e}_{i}B)^{\circ }+(\mathbf{e}_{j}B)^{\circ }+(\mathbf{e}%
_{k}B)^{\circ } \\
&\subset &(\mathbf{e}_{i}B+\mathbf{e}_{j}B+\mathbf{e}_{k}B)^{\circ }\subset
B^{\circ }.
\end{eqnarray*}

So, $B^{\circ }$ is balanced when $0\in B^{\circ }$.
\end{proof}
\end{theorem}

\begin{definition}
Let $B$ be a subset of a $H_{2}-$module $X.$ Then $B$ is called a $H-$%
absorbing set if for each $x\in X,$ there exists $\epsilon \succ 0$ such
that $\lambda x\in B$ whenever $0\preceq \lambda \preceq \epsilon .$
\end{definition}

It is clear that $H_{2}-$absorbing set always contains the origin.

\begin{theorem}
Let $B$ be a $H_{2}-$absorbing set in a $H_{2}-$module $X.$ Then for $i\in
\{1,2,3,4\},$ $\mathbf{e}_{i}B$ are absorbing sets in $%
\mathbb{R}
-$vector spaces $\mathbf{e}_{i}X.$

\begin{proof}
Let $x\in \mathbf{e}_{i}X.$ Then there exists $x^{\prime }\in X$ such that $%
x=\mathbf{e}_{i}x^{\prime }.$ Since $B$ is $H_{2}-$absorbing, there exists $%
\epsilon \succ 0$ such that $\lambda x^{\prime }\in B$ whenever $0\preceq
\lambda \preceq \epsilon .$

Let $\epsilon =\sum\limits_{i=1}^{4}\epsilon _{i}\mathbf{e}_{i}$ and $%
\lambda =\sum\limits_{i=1}^{4}\lambda _{i}\mathbf{e}_{i}.$ Then we have $%
\epsilon _{i}>0$ and $0\leq \lambda _{i}\leq \epsilon _{i}.$

Also, $\lambda _{i}x=\lambda _{i}\mathbf{e}_{i}x^{\prime }=\mathbf{e}%
_{i}\lambda _{i}x^{\prime }\in \mathbf{e}_{i}B.$

This proves that $\mathbf{e}_{i}B$ are absorbing sets in $%
\mathbb{R}
-$vector spaces $\mathbf{e}_{i}X.$
\end{proof}
\end{theorem}

\begin{remark}
We have already seen that if $B$ is $H_{2}-$balanced set, then for $i\in
\{1,2,3,4\},$ $\mathbf{e}_{i}B\subset B$ and if $B$ is $H_{2}-$convex set
containing $0,$ then $\mathbf{e}_{i}B\subset B$ . But if $B$ is $H_{2}-$%
absorbing set then $\mathbf{e}_{i}B\subset B$ does not hold.
\end{remark}

We have the following example:

\begin{example}
Let $X=H_{2}$ and $B=\{\xi \in H_{2}:\left\vert \xi \right\vert \prec \frac{1%
}{2}\}\cup \{1\}.$ Then $B$ is a $H_{2}-$absorbing of $H_{2}.$ Now $1\in B, $
so $\mathbf{e}_{i}\in \mathbf{e}_{i}B$ for $i\in \{1,2,3,4\}.$ But $\mathbf{e%
}_{i}\notin B.$ So $\mathbf{e}_{i}B\nsubseteq B.$
\end{example}

\begin{theorem}
Let $(X,\tau )$ be a topological $H_{2}-$module. Then the following
statements hold:

$(i)$ Each neighbourhood of $0$ in $X$ is $H_{2}-$absorbing.

$(ii)$ Each neighbourhood of $0$ in $X$ contains a $H_{2}-$balanced
neighbourhood of $0.$

\begin{proof}
$(i)$ Let $U\subset X$ be a neighbourhood of $0.$ Let $x\in X.$ Since scalar
multiplication is continuous and $.(0,x)=0,$ there exists a neighbourhood $V$
of $x$ and $\epsilon \succ 0$ such that whenever $\left\vert \gamma
\right\vert \prec \epsilon $ we have $\gamma V\subset U.$ In particular for $%
\gamma $ satisfying $0\preceq \gamma \preceq \frac{\epsilon }{2},$ we have $%
\gamma x\in U.$ This shows that $U$ is $H_{2}-$absorbing.

$(ii)$ Let $U\subset X$ be a neighbourhood of $0.$ Since scalar
multiplication is continuous and $.(0,0)=0,$ there exists a neighbourhood $V$
of $0$ and $\epsilon \succ 0$ such that whenever $\left\vert \gamma
\right\vert \prec \epsilon ,$ we have $\gamma V\subset U.$ Let $%
M=\bigcup\limits_{\left\vert \gamma \right\vert \prec \epsilon }\gamma V.$
Then $M$ is a neighbourhood of $0$ and $M\subset U.$ To show that $M$ is $%
H_{2}-$balanced, let $x\in M$ and $\left\vert \lambda \right\vert \preceq 1.$
Then $x=\gamma y,$ for some $y\in V.$ Since $\left\vert \lambda \gamma
\right\vert =\left\vert \lambda \right\vert \left\vert \gamma \right\vert
\prec \epsilon ,$ it follows that $\lambda x=\lambda \gamma y\in M.$\ 
\end{proof}
\end{theorem}

\section{\label{S4}Bihyperbolic-valued seminorm}

Hyperbolic-valued seminorm have been studied in $\cite{Lu14}$ and its
properties have been studied in $\cite{16}$.

In this section first we introduce bihyperbolic-valued seminorm and
investigate some properties of bihyperbolic-valued seminorm in topological
bihyperbolic modules.

\begin{definition}
Let $X$ be a $H_{2}-$module. Then a function $p:X\longrightarrow H_{2}^{+}$
is said to be a bihyperbolic-valued (or $H_{2}-$valued) seminorm if for any $%
x,y\in X$ and $\lambda \in H_{2},$ the following properties hold:

$(i)$ $p(\lambda x)=\left\vert \lambda \right\vert p(x),$

$(ii)$ $p(x+y)\preceq p(x)+p(y).$
\end{definition}

\begin{theorem}
\label{T12}Let $p$ be a $H_{2}-$valued seminorm on a $H_{2}-$module $X.$
Then for any $x,y\in X,$ the following properties hold:

$(i)$ $p(0)=0.$

$(ii)$ $\left\vert p(x)-p(y)\right\vert \preceq p(x-y).$

$(iii)$ $p(x)\succeq 0.$

$(iv)$ $\{x:p(x)=0\}$ is $H_{2}-$submodule of $X.$

\begin{proof}
The proof is similar to ($\cite{16},$ Theorem $3.2$).
\end{proof}
\end{theorem}

\begin{remark}
Every $H_{2}-$valued norm on a $H_{2}-$module is a $H_{2}-$valued seminorm.
However, the converse is not true in general.
\end{remark}

Here is an example:

\begin{example}
Define a function $p:H_{2}\longrightarrow H_{2}$ by 
\[
p(x)=\left\vert x_{1}\right\vert \mathbf{e}_{1},\text{ for each }%
x=\sum\limits_{i=1}^{4}x_{i}\mathbf{e}_{i}\in H_{2}. 
\]

Then clearly $p$ is a $H_{2}-$valued seminorm on $H_{2}.$ Now $\mathbf{e}%
_{2}\in H_{2}$ and $\mathbf{e}_{2}\neq 0,$ but $p(\mathbf{e}_{2})=0.$

So, $p$ is not a $H_{2}-$valued norm on $H_{2}.$
\end{example}

\begin{theorem}
Let $p$ be a $H_{2}-$valued seminorm on a topological $H_{2}-$module $X.$
Denote the sets $\{x\in X:p(x)\prec 1\}$ and $\{x\in X:p(x)\preceq 1\}$ by $%
A $ and $C$ respectively. Then the following statements are equivalent:

$(i)$ $p$ is continuous.

$(ii)$ $A$ is open.

$(iii)$ $0\in A^{\circ }.$

$(iv)$ $0\in C^{\circ }.$

$(v)$ $p$ is continuous at $0.$

$(vi)$ there exists a continuous $H_{2}-$valued seminorm $q$ on $X$ such
that $p\preceq q.$

\begin{proof}
The proof is similar to ($\cite{16},$ Theorem $3.5$).
\end{proof}
\end{theorem}

\begin{theorem}
\label{T14}Let $X$ be a $H_{2}-$module and $p$ be a $H_{2}-$valued seminorm
on $X.$ Then $\{x\in X:p(x)\prec 1\}$ and $\{x\in X:p(x)\preceq 1\}$ are $%
H_{2}-$convex, $H_{2}-$balanced and $H_{2}-$absorbing on $X.$

\begin{proof}
The proof is similar to ($\cite{16},$ Theorem $3.6$).
\end{proof}
\end{theorem}

\section{\label{S5}Bihyperbolic-valued Minkowski Functionals}

Hyperbolic-valued Minkowski functionals in hyperbolic modules have been
studied in $\cite{Lu14}$ and hyperbolic-valued Minkowski functionals in
bicomplex modules have been studied in $\cite{16}$.

In this section, we define bihyperbolic-valued Minkowski functionals in
bihyperbolic modules and it has been shown that a $H_{2}-$valued Minkowski
functionals of a $H_{2}-$balanced, $H_{2}-$convex and $H_{2}-$absorbing set
turns out to be a $H_{2}-$valued seminorm.

\begin{definition}
Let $B$ be a $H_{2}-$convex, $H_{2}-$absorbing subset of a $H_{2}-$module $%
X. $ Then the mapping $q_{B}:X\longrightarrow $ $H_{2}^{+}$ defined by%
\[
q_{B}(x)=\inf_{H_{2}}\{\alpha \succ 0:x\in \alpha B\},\text{ \ for each }%
x\in X 
\]

is called bihyperbolic-valued gauge or bihyperbolic-valued Minkowski
functional on $B.$
\end{definition}

Since $B$ is $H_{2}-$convex, then by Theorem \ref{T4}, we have $%
B=\sum\limits_{i=1}^{4}\mathbf{e}_{i}B.$ Then for $x=\sum\limits_{i=1}^{4}%
\mathbf{e}_{i}x_{i},$ $\alpha =\sum\limits_{i=1}^{4}\mathbf{e}_{i}\alpha
_{i},$ $q_{B}$ can be written as%
\[
q_{B}(x)=\sum\limits_{i=1}^{4}\mathbf{e}_{i}q_{Bi}(x), 
\]

where 
\[
q_{Bi}(x)=\inf \{\alpha _{i}>0:x_{i}\in \alpha _{i}\mathbf{e}_{i}B\}. 
\]

\begin{theorem}
\label{T15}Let $B$ be a $H_{2}-$convex, $H_{2}-$balanced, $H_{2}-$absorbing
subset of a $H_{2}-$module $X.$ Then the $H_{2}-$valued gauge $q_{B}$ is a $%
H_{2}-$valued seminorm on $X.$

\begin{proof}
Let $x,y\in X$ such that $q_{B}(x)=\alpha $ and $q_{B}(y)=\gamma .$ Then for
any $\epsilon \succ 0,$ we have $x\in (\alpha +\epsilon )B$ and $y\in
(\gamma +\epsilon )B.$ Therefore we can find $u,v\in B$ such that $x=(\alpha
+\epsilon )u$ and $y=(\gamma +\epsilon )v.$ Observe that $0\prec \frac{%
(\alpha +\epsilon )}{(\alpha +\gamma +2\epsilon )}\prec 1$ and $0\prec \frac{%
(\gamma +\epsilon )}{(\alpha +\gamma +2\epsilon )}\prec 1.$ Therefore, by $%
H_{2}-$convexity of $B,$ we have%
\[
\frac{(\alpha +\epsilon )u+(\gamma +\epsilon )v}{(\alpha +\gamma +2\epsilon )%
}\in B\Longrightarrow (\alpha +\epsilon )u+(\gamma +\epsilon )v\in (\alpha
+\gamma +2\epsilon )B 
\]

which implies that 
\[
x+y\in (\alpha +\gamma +2\epsilon )B. 
\]

Letting $\epsilon \longrightarrow 0,$ we obtain%
\[
q_{B}(x+y)\preceq \alpha +\gamma =q_{B}(x)+q_{B}(y). 
\]%
\newline

We now show that $q_{B}(\lambda x)=\left\vert \lambda \right\vert $ $%
q_{B}(x) $ for each $x\in X$ and $\lambda \in H_{2}.$

Clearly, $q_{B}(0)=0.$ So we assume that $\lambda \in H_{2}\backslash NC.$
Since $B$ is $H_{2}-$balanced, by Theorem \ref{T1}, we have%
\begin{eqnarray*}
q_{B}(\lambda x) &=&\inf_{H_{2}}\{\alpha \succ 0:\lambda x\in \alpha B\} \\
&=&\inf_{H_{2}}\left\{ \alpha \succ 0:x\in \alpha \left( \frac{1}{\lambda }%
B\right) \right\} \\
&=&\inf_{H_{2}}\left\{ \alpha \succ 0:x\in \alpha \left( \frac{1}{\left\vert
\lambda \right\vert }B\right) \right\} \\
&=&\left\vert \lambda \right\vert \inf_{H_{2}}\left\{ \frac{\alpha }{%
\left\vert \lambda \right\vert }\succ 0:x\in \frac{\alpha }{\left\vert
\lambda \right\vert }B\right\} \\
&=&\left\vert \lambda \right\vert q_{B}(x).
\end{eqnarray*}

Now suppose $\lambda \in NC^{\ast }$ such that $\lambda =\lambda _{i}\mathbf{%
e}_{i}$ for $i\in \{1,2,3,4\}.$ Since $B$ is $H_{2}-$balanced set in $X,$ it
follows that $\mathbf{e}_{i}B$ is balanced set in $%
\mathbb{R}
-$vector space $\mathbf{e}_{i}X.$ Hence%
\begin{eqnarray*}
q_{B}(\lambda x) &=&\sum\limits_{l=1}^{4}\mathbf{e}_{l}q_{Bl}(\lambda x) \\
&=&\mathbf{e}_{i}\inf \{\alpha _{i}>0:\lambda _{i}x_{i}\in \alpha _{i}%
\mathbf{e}_{i}B\} \\
&=&\mathbf{e}_{i}\inf \left\{ \alpha _{i}>0:x_{i}\in \alpha _{i}\left( \frac{%
1}{\lambda _{i}}\mathbf{e}_{i}B\right) \right\} \\
&=&\mathbf{e}_{i}\inf \left\{ \alpha _{i}>0:x_{i}\in \alpha _{i}\left( \frac{%
1}{\left\vert \lambda _{i}\right\vert }\mathbf{e}_{i}B\right) \right\} \\
&=&\left\vert \lambda _{i}\right\vert \mathbf{e}_{i}\inf \left\{ \frac{%
\alpha _{i}}{\left\vert \lambda _{i}\right\vert }>0:x_{i}\in \frac{\alpha
_{i}}{\left\vert \lambda _{i}\right\vert }\mathbf{e}_{i}B\right\} \\
&=&\left\vert \lambda _{i}\right\vert \mathbf{e}_{i}q_{Bi}(x)=\left\vert
\lambda \right\vert q_{B}(x).
\end{eqnarray*}

Now suppose $\lambda \in NC^{\ast }$ such that $\lambda =\lambda _{i}\mathbf{%
e}_{i}+\lambda _{j}\mathbf{e}_{j}$ for $i\neq j,$ $i,j\in \{1,2,3,4\}.$
Since $B$ is $H_{2}-$balanced set in $X,$ it follows that $\mathbf{e}_{i}B$,$%
\mathbf{e}_{j}B$ are balanced set in $%
\mathbb{R}
-$vector spaces $\mathbf{e}_{i}X,\mathbf{e}_{j}X$ respectively. Hence%
\begin{eqnarray*}
q_{B}(\lambda x) &=&\sum\limits_{l=1}^{4}\mathbf{e}_{l}q_{Bl}(\lambda x) \\
&=&\mathbf{e}_{i}\inf \{\alpha _{i}>0:\lambda _{i}x_{i}\in \alpha _{i}%
\mathbf{e}_{i}B\}+\mathbf{e}_{j}\inf \{\alpha _{j}>0:\lambda _{j}x_{j}\in
\alpha _{j}\mathbf{e}_{j}B\} \\
&=&\mathbf{e}_{i}\inf \left\{ \alpha _{i}>0:x_{i}\in \alpha _{i}\left( \frac{%
1}{\lambda _{i}}\mathbf{e}_{i}B\right) \right\} +\mathbf{e}_{j}\inf \left\{
\alpha _{j}>0:x_{j}\in \alpha _{j}\left( \frac{1}{\lambda _{j}}\mathbf{e}%
_{j}B\right) \right\} \\
&=&\mathbf{e}_{i}\inf \left\{ \alpha _{i}>0:x_{i}\in \alpha _{i}\left( \frac{%
1}{\left\vert \lambda _{i}\right\vert }\mathbf{e}_{i}B\right) \right\} +%
\mathbf{e}_{j}\inf \left\{ \alpha _{j}>0:x_{j}\in \alpha _{j}\left( \frac{1}{%
\left\vert \lambda _{j}\right\vert }\mathbf{e}_{j}B\right) \right\} \\
&=&\left\vert \lambda _{i}\right\vert \mathbf{e}_{i}\inf \left\{ \frac{%
\alpha _{i}}{\left\vert \lambda _{i}\right\vert }>0:x_{i}\in \frac{\alpha
_{i}}{\left\vert \lambda _{i}\right\vert }\mathbf{e}_{i}B\right\}
+\left\vert \lambda _{j}\right\vert \mathbf{e}_{j}\inf \left\{ \frac{\alpha
_{j}}{\left\vert \lambda _{j}\right\vert }>0:x_{j}\in \frac{\alpha _{j}}{%
\left\vert \lambda _{j}\right\vert }\mathbf{e}_{j}B\right\} \\
&=&\left\vert \lambda _{i}\right\vert \mathbf{e}_{i}q_{Bi}(x)+\left\vert
\lambda _{j}\right\vert \mathbf{e}_{j}q_{Bj}(x)=\left\vert \lambda
\right\vert q_{B}(x).
\end{eqnarray*}

Now suppose $\lambda \in NC^{\ast }$ such that $\lambda =\lambda _{i}\mathbf{%
e}_{i}+\lambda _{j}\mathbf{e}_{j}+\lambda _{k}\mathbf{e}_{k}$ for $i\neq
j\neq k,$ $i,j,k\in \{1,2,3,4\}.$ Since $B$ is $H_{2}-$balanced set in $X,$
it follows that $\mathbf{e}_{i}B$,$\mathbf{e}_{j}B,\mathbf{e}_{k}B$ are
balanced set in $%
\mathbb{R}
-$vector spaces $\mathbf{e}_{i}X,\mathbf{e}_{j}X$,$\mathbf{e}_{k}X$
respectively. Hence by similar technique used as above, we can prove in this
case also%
\[
q_{B}(\lambda x)=\left\vert \lambda \right\vert q_{B}(x). 
\]

This completes the proof.
\end{proof}
\end{theorem}

\begin{definition}
Let $(X,\tau )$ be a topological $H_{2}-$module. Then a subset $B\subset X$
is said to be bounded if for each neighbourhood $U$ of $0,$ there exists $%
\lambda \succ 0$ such that $B\subset \lambda U.$
\end{definition}

\begin{corollary}
Let $B$ be a bounded $H_{2}-$convex, $H_{2}-$balanced, $H_{2}-$absorbing
subset of a topological $H_{2}-$module $(X,\tau ).$ Then $q_{B}$ is $H_{2}-$%
valued norm on $X.$
\end{corollary}

The next result follows from Theorems \ref{T14} and \ref{T15}.

\begin{theorem}
Let $B$ be a $H_{2}-$convex, $H_{2}-$balanced, $H_{2}-$absorbing subset of a 
$H_{2}-$module $X$ and $q_{B}$ be the $H_{2}-$valued gauge on $B.$ Then both
the subsets $\{x\in X:$ $q_{B}(x)\prec 1\}$ and $\{x\in X:$ $q_{B}(x)\preceq
1\}$ of $X$ are $H_{2}-$convex, $H_{2}-$balanced and $H_{2}-$absorbing.
\end{theorem}

\begin{theorem}
Let $(X,\tau )$ be a topological $H_{2}-$module, $B$ be a $H_{2}-$convex, $%
H_{2}-$balanced, $H_{2}-$absorbing subset of $X$ and $q_{B}$ be the $H_{2}-$%
valued gauge on $B.$ Let us denote $\{x\in X:$ $q_{B}(x)\prec 1\}$ and $%
\{x\in X:$ $q_{B}(x)\preceq 1\}$ by $A_{B}$ and $C_{B}$ respectively. Then,
the following statements hold:

$(i)$ $B^{\circ }\subset A_{B}\subset B\subset C_{B}\subset \overline{B}.$

$(ii)$ If $B$ is open, then $B=A_{B}.$

$(iii)$ If $B$ is closed, then $B=C_{B}.$

$(iv)$ If $q_{B}$ is continuous, then $B^{\circ }=A_{B}.$

\begin{proof}
The proof is similar to ($\cite{16},$ Theorem $4.6$).
\end{proof}
\end{theorem}

\section{\label{S6}Locally Bihyperbolic Convex Modules}

In this section, we introduce the bihyperbolic version of locally convex
topological spaces, bihyperbolic metrizable and bihyperbolic normable
locally bihyperbolic convex modules.

\begin{definition}
Let $X$ be a $H_{2}-$module and $\mathcal{P}$ be a family of $H_{2}-$valued
seminorms on $X$. Then the family $\mathcal{P}$ is said to be separated if
for each $x\neq 0$, there exists $p\in \mathcal{P}$ such that $p(x)\neq 0.$
\end{definition}

We define a topology on a $H_{2}-$module $X$ determined by the family $%
\mathcal{P}$ of $H_{2}-$valued seminorms on $X$ as follows:

For $x\in X$, $\epsilon \succ 0$ and $p\in \mathcal{P}$, we set%
\[
U(x,\epsilon ,p)=\{y\in x:p(y-x)\prec \epsilon \}, 
\]

and for $x\in X$, $\epsilon \succ 0$ and $p_{1},p_{2},...,p_{n}\in \mathcal{P%
},$ set 
\[
U(x,\epsilon ,p_{1},p_{2},...,p_{n})=\{y\in x:p_{i}(y-x)\prec \epsilon
,p_{2}(y-x)\prec \epsilon ,...,p_{n}(y-x)\prec \epsilon \}. 
\]

Let $\mathcal{U}_{\mathcal{P}}(x)=\{U(x,\epsilon
,p_{1},p_{2},...,p_{n}):\epsilon \succ 0,$ $p_{1},p_{2},...,p_{n}\in 
\mathcal{P}$ and $n\in 
\mathbb{N}
\}.$

Then $\mathcal{U}_{\mathcal{P}}=\{\mathcal{U}_{\mathcal{P}}(x):x\in
X\}=\bigcup\limits_{x\in X}\mathcal{U}_{\mathcal{P}}(x)$ forms a base for a
topology $\tau _{\mathcal{P}}$ on $X,$ called topology generated by the
family $\mathcal{P}$.

\begin{theorem}
Let $\mathcal{P}$ be a separated family of $H_{2}-$valued seminorms on $X.$
Then $(X,\tau _{\mathcal{P}})$ is a topological $H_{2}-$module.

\begin{proof}
The proof is similar to ($\cite{16},$ Theorem $5.2$).
\end{proof}
\end{theorem}

\begin{lemma}
Let $X$ is a topological $H_{2}-$module and $\mathcal{P=\{}p_{n}\}_{n\in 
\mathbb{N}
}$ be a family of $H_{2}-$valued seminorms on $X.$ For each $m\in 
\mathbb{N}
,$ define $q_{m}:X\longrightarrow H_{2}$ by 
\[
q_{m}(x)=\sup \{p_{1}(x),p_{2}(x),...,p_{m}(x)\}\text{ for each }x\in X. 
\]

Then,$\mathcal{Q=\{}q_{m}\}_{m\in 
\mathbb{N}
}$ is a family of $H_{2}-$valued seminorms on $X$ such that the following
hold:

$(i)$ $\mathcal{Q}$ is separated if $\mathcal{P}$ is so.

$(ii)$ $q_{m}\preceq q_{m+1}$ for each $m\in 
\mathbb{N}
.$

$(iii)$ $(X,\tau _{\mathcal{P}})$ and $(X,\tau _{\mathcal{Q}})$ are
topologically isomorphic.

\begin{proof}
The proof is similar to ($\cite{La}$, Lemma $2.5.1$).
\end{proof}
\end{lemma}

\begin{definition}
A topological $H_{2}-$module $(X,\tau )$ is said to be locally bihyperbolic
convex (or $H_{2}-$convex) module if it has a neighbourhood base at $0$ of $%
H_{2}-$convex sets.
\end{definition}

\begin{theorem}
A topological $H_{2}-$module $(X,\tau )$ is a locally $H_{2}-$convex module
if and only if its topology is generated by a separated family $\mathcal{P}$
of $H_{2}-$valued seminorms on $X.$

\begin{proof}
The proof is similar to ($\cite{16},$ Theorem $5.5$).
\end{proof}
\end{theorem}

\begin{definition}
Let $d_{H_{2}}:$ $X\times X$ $\rightarrow H_{2}$ be a function such that for
any $x,y,z\in X,$ the following properties hold:

$(i)$ $d_{H_{2}}(x,y)\succeq $ $0$ and $d_{H_{2}}(x,y)=0$ if and only if $%
x=y,$

$(ii)$ $d_{H_{2}}(x,y)=d_{H_{2}}(y,x),$

$(iii)$ $d_{H_{2}}(x,z)\preceq d_{H_{2}}(x,y)+d_{H_{2}}(y,z).$

Then $d_{H_{2}}$ is called a bihyperbolic-valued (or $H_{2}-$valued) metric
on $X$ and the pair $(X,d_{H_{2}})$ is called a bihyperbolic metric (or $%
H_{2}-$metric) space.
\end{definition}

The following result is easy to prove

\begin{lemma}
Every $H_{2}-$metric space is first countable.
\end{lemma}

\begin{definition}
A topological $H_{2}-$module $X$ is said to be bihyperbolic metrizable (or $%
H_{2}-$metrizable) if the topology on $X$ is generated by a $H_{2}-$valued
metric on $X.$
\end{definition}

\begin{definition}
A topological $H_{2}-$module $X$ is said to be bihyperbolic normable (or $%
H_{2}-$normable) if the topology on $X$ is generated by a $H_{2}-$valued
norm on $X.$
\end{definition}

\begin{lemma}
Let $\mathcal{P}=\{p_{n}\}$ be a countable separated family of $H_{2}-$%
valued seminorms on a topological $H_{2}-$module $(X,\tau )$ such that $%
p_{n}\preceq p_{n+1}$ for each $n\in 
\mathbb{N}
$. Define a function $d:$ $X\times X$ $\rightarrow H_{2}$ by 
\[
d(x,y)=\sum\limits_{n=1}^{\infty }2^{-n}\frac{p_{n}(x-y)}{1+p_{n}(x-y)},%
\text{ for each }x,y\in X. 
\]

Then, $d$ is a translation invariant $H_{2}-$valued metric on $X$ and the
topology on $X$ generated by $d$ is the topology generated by the family $%
\mathcal{P}$.

\begin{proof}
The proof is similar to ($\cite{La}$, Theorem $2.5.1$).
\end{proof}
\end{lemma}

\begin{theorem}
A locally $H_{2}-$convex module $(X,\tau )$ is $H_{2}-$metrizable if and
only if its topology is generated by a countable separated family $\mathcal{P%
}$ of $H_{2}-$valued seminorms on $X$.

\begin{proof}
The proof is similar to ($\cite{16},$ Theorem $5.11$).
\end{proof}
\end{theorem}

\begin{theorem}
A topological $H_{2}-$module $(X,\tau )$ is $H_{2}-$normable if and only if
it contains a bounded $H_{2}-$convex neighbourhood of $0.$

\begin{proof}
The proof is similar to ($\cite{16},$ Theorem $5.12$).
\end{proof}
\end{theorem}

\end{document}